\documentclass[10pt,twoside]{article}
\usepackage{graphicx}
\usepackage{amsmath}
\usepackage{Latex-document}

\input amssym.def
\input amssym.tex

\def\bsh{\backslash}
\def\bu{\bullet}
\def\es{\emptyset}
\def\fl{\forall \, }
\def\hra{\hookrightarrow}
\def\ify{\infty}
\def\longra{\longrightarrow}
\def\mpo{\mapsto}
\def\op{\oplus}
\def\ot{\otimes}
\def\ov{\overline}
\def\ra{\rightarrow}
\def\rra{\rightrightarrows}
\def\sbs{\subset}
\def\sbsq{\subseteq}
\def\thra{\twoheadrightarrow}
\def\ts{\times}
\def\ud{\underline}

\def\lhra{\lhook\joinrel\relbar\joinrel\rightarrow}
\def\lhla{\leftarrow\joinrel\relbar\joinrel\rhook}

\font\tenbb=msbm10
\font\sevenbb=msbm7
\font\fivebb=msbm5
\newfam\bbfam
\textfont\bbfam=\tenbb
\scriptfont\bbfam=\sevenbb
\scriptscriptfont\bbfam=\fivebb
\def\bb{\fam\bbfam}

\def\Ab{{\bb A}}
\def\Cb{{\bb C}}
\def\Fb{{\bb F}}
\def\Gb{{\bb G}}
\def\Nb{{\bb N}}
\def\Qb{{\bb Q}}
\def\Rb{{\bb R}}
\def\Zb{{\bb Z}}

\def\Aut{\mathop{\rm Aut}\nolimits}
\def\Cht{\mathop{\rm Cht}\nolimits}

\def\ess{\mathop{\rm ess}\nolimits}
\def\Frob{\mathop{\rm Frob}\nolimits}
\def\GL{\mathop{\rm GL}\nolimits}
\def\Id{\mathop{\rm Id}\nolimits}
\def\Ker{\mathop{\rm Ker}\nolimits}
\def\Lef{\mathop{\rm Lef}\nolimits}
\def\Spec{\mathop{\rm Spec}\nolimits}
\def\Tr{\mathop{\rm Tr}\nolimits}

\def\d{\delta}
\def\g{\gamma}
\def\G{\Gamma}
\def\k{\kappa}
\def\lb{\lambda}
\def\Om{\Omega}
\def\s{\sigma}
\def\vp{\varphi}

\def\Ec{{\cal E}}
\def\Hc{{\cal H}}

\def\Oc{{\cal O}}
\def\Vc{{\cal V}}

\def\build#1_#2^#3{\mathrel{
\mathop{\kern 0pt#1}\limits_{#2}^{#3}}}

\def\un{{\rm 1\mkern-4mu I}}

\def\varb{\varsubsetneq}

\def\thebibliography#1{\section*{Bibliographie}
\list{[\arabic{enumi}]}{\settowidth \labelwidth{[#1]} \leftmargin
\labelwidth \advance \leftmargin \labelsep \usecounter{enumi}}
\def\newblock{\hskip .11em plus .33em minus .07em} \sloppy
\clubpenalty 4000 \widowpenalty 4000 \sfcode`\.=1000 \relax}

\font\goth=eufm10

\def\volumeno{I}

\title{\bf  Chtoucas de Drinfeld, Formule des \vskip -2mm Traces d'Arthur-Selberg et\vskip -2mm Correspondance de Langlands\vskip 6mm}

\author{Laurent Lafforgue\vspace*{-0.5cm}\thanks{Institut des Hautes Etudes
Scientifiques, Le Bois Marie, 35 route de Chartres, F-91440 \break Bures-sur-Yvette,
France. Courriel : laurent@ihes.fr}}
\date{\vspace{-8mm}}

\begin{document}

\maketitle

\thispagestyle{first} \setcounter{page}{383}

\begin{quotation}
\centerline{\bf R\'esum\'e}

\vskip 3mm

Cet expos\'e est consacr\'e \`a la preuve de la correspondance de
Lan-glands pour les groupes $\GL_r$ sur les corps de fonctions.

\vskip 4mm

\noindent {\bf Code mati\`ere (AMS 2000):} 11F, 11F52, 11F60,
11F66, 11F70, 11F72, 11F80, 11R39, 14G35, 14H60, 22E55.

\noindent {\bf Mots-Cl\'es:} Chtoucas, Vari\'et\'es modulaires de
Drinfeld, Modules des fibr\'es sur les courbes, Correspondance de
Langlands, Corps de fonctions, Repr\'esentations galoi\-siennes,
Repr\'esentations automorphes, Fonctions L, Cohomologie
$\ell$-adique, Correspondances de Hecke, Formule des traces
d'Arthur-Selberg, Formule des points fixes de
Grothendieck-Lefschetz.
\end{quotation}

\vskip 10mm

\section*{Introduction} \label{section 0}\setzero
\vskip-5mm \hspace{5mm }

La g\'eom\'etrie alg\'ebrique sur un corps $F$ consiste \`a \'etudier les
vari\'et\'es alg\'ebriques sur $F$ (c'est-\`a-dire principalement les objets
g\'eom\'etriques d\'efinis dans des espaces projectifs par des syst\`emes
d'\'equations polynomiales \`a coefficients dans $F$) et les morphismes ou plus
g\'en\'eralement les correspondances qui les relient.

Si $\ell$ est un nombre premier diff\'erent de la caract\'eristique de $F$,
Grothen\-dieck a associ\'e \`a toute vari\'et\'e alg\'ebrique $V$ sur $F$ ses
espaces de cohomologie $\ell$-adique $H_c^{\nu} (V \ot \ov F , \Qb_{\ell})$, $0 \leq
\nu \leq 2 \dim V$, qui sont des repr\'esentations du groupe de Galois $G_F$ de $F$
continues et de dimension finie sur $\Qb_{\ell}$. Quand le corps $F$ est de type
fini sur $\Qb$ ou $\Fb_p = \Zb / p\Zb$, la th\'eorie conjecturale des motifs de
Grothendieck et les conjectures de Tate pr\'evoient qu'il devrait \^etre possible
de remonter des repr\'esentations $\ell$-adiques de $G_F$ aux vari\'et\'es
projectives lisses sur $F$ et \`a leurs correspondances si bien que la
connaissance de celles-ci serait essentiellement \'equivalente \`a celle des
repr\'esentations $\ell$-adiques irr\'eductibles de $G_F$.

Cette perspective pose de fa\c con cruciale la question de dresser la liste de
toutes les repr\'esentations $\ell$-adiques irr\'eductibles de $G_F$. Quand $F$
est un ``corps global'', c'est-\`a-dire une extension finie de $\Qb$ ou bien le
corps $F(X)$ des fonctions rationnelles d'une courbe $X$ projective lisse sur un
corps fini, Langlands a conjectur\'e que les repr\'esentations irr\'eductibles de
dimension $r$ de $G_F$ sont param\'etr\'ees naturellement par les
repr\'esentations automorphes cuspidales de $\GL_r$ sur $F$.

L'expos\'e est consacr\'e \`a la preuve de la correspondance de Langlands dans le
cas des corps de fonctions $F=F(X)$. Un r\^ole central y est jou\'e par la formule
des traces d'Arthur-Selberg que la th\'eorie des chtoucas de Drinfeld permet
d'interpr\'eter g\'eom\'etriquement.

\section{L'\'enonc\'e de la correspondance de Langlands} \label{section 1}\setzero
\vskip-5mm \hspace{5mm }

On consid\`ere donc une courbe $X$ projective, lisse et g\'eom\'etriquement
connexe sur un corps $\Fb_q$ \`a $q$ \'el\'ements. On note $F=F(X)$ le corps des
fonctions rationnelles sur $X$ et $\vert X \vert$ l'ensemble des points ferm\'es
de $X$. Pour $x \in \vert X \vert$, on note $\k_x$ son corps r\'esiduel et $\deg
(x)$ la dimension de $\k_x$ sur $\Fb_q$.

\subsection{\boldmath Repr\'esentations galoisiennes $\ell$-adiques}\label{ssec1a}\setzero
\vskip-5mm \hspace{5mm }

On note $G_F$ le groupe de Galois de $F$, c'est-\`a-dire le groupe des
automorphismes d'une cl\^oture s\'eparable $\ov F$ de $F$. C'est un groupe
profini.

Soit $\ell$ un nombre premier diff\'erent de la caract\'eristique de $\Fb_q$.

Une repr\'esentation $\ell$-adique $\s$ de $G_F$ est une repr\'esentation
continue et de dimension finie sur $\ov{\Qb}_{\ell}$ qui est d\'efinie sur une
extension finie de $\Qb_{\ell}$ et provient d'un faisceau $\ell$-adique lisse (ou
``syst\`eme local'') sur un ouvert non vide de $X$. On dit que $\s$ est non
ramifi\'ee en un point $x \in \vert X \vert$ si $x$ est dans le plus grand ouvert
o\`u $\s$ est un syst\`eme local. Dans ce cas, la fibre $\s_x$ de $\s$ en $x$ est
une repr\'esentation du groupe de Galois de $\k_x$ lequel est engendr\'e par
l'automorphisme $\Frob_x$ d'\'el\'evation \`a la puissance $q^{\deg (x)}$ et on
peut poser
$$
{\rm L}_x (\s , T) = \textstyle{\det_{\s_x}} (1-T \cdot \Frob_x^{-1})^{-1} \, .
$$

Pour tout entier $r \geq 1$, on note $\{ \s \}_r$ l'ensemble des classes
d'isomorphie de repr\'esentations $\ell$-adiques irr\'eductibles de dimension
$r$ de $G_F$ dont le d\'eterminant est d'ordre fini.

\subsection{Groupes ad\'eliques et alg\`ebres de Hecke}\label{ssec1b}\setzero
\vskip-5mm \hspace{5mm }

Tout point $x \in \vert X \vert$ induit sur $F$ la valuation
$$
\deg_x : F^{\ts} \mpo \Zb \qquad f \mpo \deg_x (f) = \ \hbox{l'ordre d'annulation de
$f$ en $x$.}
$$
On note $F_x$ le corps compl\'et\'e de $F$ en $x$~; il contient comme sous-anneau
d'entiers $O_x = \{ f_x \in F_x \mid \deg_x (f_x) \geq 0 \}$ qui est compact.

L'anneau $\Ab$ des ad\`eles de $F$ est le produit ``restreint''
$\build{\prod\!\!\!\!\!\!\!\coprod}_{x \, \in \, \vert X \vert}^{} F_x$ des
familles $f_x \in F_x$, $x \in \vert X \vert$, telles que $f_x \in O_x$ pour
presque tout $x$ (c'est-\`a-dire tout $x$ sauf un nombre fini). Il contient comme
sous-anneau compact des entiers $O_{\Ab} = \build\prod_{x \, \in \, \vert X \vert}^{}
O_x = \ \build\varprojlim_{N}^{} \Oc_N$, o\`u $N$ d\'ecrit l'ensemble des ``niveaux''
c'est-\`a-dire des sous-sch\'emas ferm\'es finis $N = \Spec \Oc_N \hra X$.

L'anneau $\Ab$ contient $F$ comme sous-groupe additif discret cocompact et le
noyau $\Ab^{\ts 0}$ de l'homomorphisme
$$
\deg : \Ab^{\ts} \ra \Zb \qquad (f_x) \mpo \sum_{x \, \in \, \vert X \vert} \deg
(x) \, \deg_x (f_x)
$$
contient $F^{\ts}$ comme sous-groupe discret cocompact.

Pour tout entier $r \geq 1$, le groupe topologique $\GL_r (F_x)$ [resp. $\GL_r
(\Ab) = \build{\prod\!\!\!\!\!\!\!\coprod}_{x \, \in \, \vert X \vert}^{} \GL_r
(F_x)$] a une mesure de Haar ${\rm d} g_x$ [resp. ${\rm d} g_{\Ab} = \prod {\rm d}
g_x$]
qui attribue le volume 1 au sous-groupe ouvert compact maximal $K_x = \GL_r (O_x)$
[resp. $K = \GL_r (O_{\Ab}) =  \prod K_x$]. L'alg\`ebre de Hecke est l'alg\`ebre
de convolution $\Hc_x^r$ [resp. $\Hc^r = \ot \, \Hc_x^r$] des fonctions localement
constantes \`a support compact sur $\GL_r (F_x)$ [resp. $\GL_r (\Ab)$]. Elle est
r\'eunion filtrante des sous-alg\`ebres $\Hc_{N,x}^r$ [resp. $\Hc_N^r = \ot \,
\Hc_{N,x}^r$] des fonctions bi-invariantes par les sous-groupes ouverts d'indice
fini $K_{N,x} = \Ker (K_x \ra \GL_r (\Oc_N))$ [resp. $K_N = \prod K_{N,x} = \Ker
(K \thra \GL_r (\Oc_N)$]. Chaque $\Hc_{N,x}^r$ [resp. $\Hc_N^r$] a une unit\'e
$\un_{N,x}$ [resp. $\un_N = \ot \, \un_{N,x}$].

On s'int\'eresse aux repr\'esentations ``lisses admissibles irr\'eductibles'' de
$\GL_r (F_x)$ [resp. $\GL_r (\Ab)$]. Ce sont les modules simples $\pi_x$ sur
$\Hc_x^r$ [resp. $\pi$ sur $\Hc^r$] qui sont r\'eunions filtrantes des $\pi_x
\cdot \un_{N,x}$ [resp. $\pi \cdot \un_N$] suppos\'es de dimension finie. Chaque
$\pi_x \cdot \un_{N,x}$ [resp. $\pi \cdot \un_N$] est un module sur $\Hc_{N,x}^r$
[resp. $\Hc_N^r$], et s'il n'est pas nul il est simple et caract\'erise $\pi_x$
[resp. $\pi$]. L'unit\'e $\un_{\es , x}$ de l'alg\`ebre de Hecke ``sph\'erique''
$\Hc_{\es , x}^r$ est la fonction caract\'eristique de $K_x = \GL_r (O_x)$;
quand $\pi_x \cdot \un_{\es , x} \ne 0$, on dit que $\pi_x$ est ``non
ramifi\'ee''. Enfin, les repr\'esentations lisses admissibles irr\'eductibles de
$\GL_r (\Ab)$ sont celles de la forme $\pi = \build\bigotimes_{x \, \in \, \vert X
\vert}^{} \pi_x$ o\`u les $\pi_x$ sont des repr\'esentations lisses admissibles
irr\'eductibles des $\GL_r (F_x)$ presque toutes non ramifi\'ees.

\smallskip

{\bf Th\'eor\`eme (Satake).} -- \it  L'alg\`ebre sph\'erique $\Hc_{\es , x}^r$ est
commutative et isomorphe \`a l'alg\`ebre des polyn\^omes sym\'etriques en $Z_1 ,
Z_1^{-1} , \ldots , Z_r , Z_r^{-1}$. \rm

\smallskip

Par cons\'equent, une repr\'esentation lisse admissible non
ramifi\'ees irr\'eductible $\pi_x$ de $\GL_r (F_x)$ est
caract\'eris\'ee par $r$ scalaires non nuls $z_1 (\pi_x) , \ldots
, z_r (\pi_x)$ (appel\'es valeurs propres de Hecke) bien d\'efinis
\`a l'ordre pr\`es ou par
$$
{\rm L} (\pi_x , T) = \prod_{i \, = \, 1}^r (1 - T \cdot z_i (\pi_x)^{-1})^{-1} \, .
$$
Et si $\pi = \ot \, \pi_x$ est une repr\'esentation lisse admissible
irr\'eductible de $\GL_r (\Ab)$, ${\rm L}_x (\pi , T) = {\rm L} (\pi_x , T)$ est
d\'efini d\'ej\`a en tout $x$ o\`u $\pi_x$ est non ramifi\'ee, donc en presque tout
$x$.

\newpage

\subsection{Repr\'esentations automorphes cuspidales et d\'ecomposi\-tion spectrale de
Langlands}\label{ssec1c}\setzero \vskip-5mm \hspace{5mm }

Pour tout entier $r \geq 1$, $\GL_r (F)$ est un sous-groupe discret de covolume
fini (mais non cocompact si $r \geq 2$) dans $\GL_r (\Ab)^0 = \Ker (\GL_r (\Ab)
\build\longra_{}^{\det} \Ab^{\ts} \build\longra_{}^{\deg} \Zb)$. La th\'eorie
automorphe consiste \`a \'etudier le quotient
$$
\GL_r (F) \bsh \GL_r (\Ab)
$$
via l'espace de ses fonctions muni de l'action de $\Hc^r$ par convolution \`a
droite.

Consid\'erons d'abord l'espace des fonctions $\vp$ sur $\GL_r (\Ab)$ qui sont
localement constantes, \`a support compact, invariantes \`a gauche par $\GL_r (F)$
et v\'erifient $\vp (ag) = \vp (g)$, $\fl g$, pour un certain $a \in \Ab^{\ts}$ de
degr\'e $\ne 0$ et la condition de ``cuspidalit\'e''
$$
\int_{N_P (F) \bsh N_P (\Ab)} \vp (n_P \, g) \cdot {\rm d} n_P = 0 \, , \quad \fl g \in
\GL_r (\Ab) \, , \quad \fl P \varb \GL_r \, ,
$$
o\`u $P$ d\'ecrit l'ensemble des sous-groupes paraboliques standard de $\GL_r$,
$N_P$ d\'esigne le radical unipotent de $P$ et ${\rm d} n_P$ une mesure de Haar sur
$N_P
(\Ab)$.

Il s'\'ecrit comme une somme directe de repr\'esentations lisses admissibles
irr\'educ\-tibles de $\GL_r (\Ab)$ appel\'ees les repr\'esentations automorphes
cuspidales. On note $\{ \pi \}_r$ leur ensemble. Chaque $\pi \in \{ \pi \}_r$ a un
caract\`ere central $\chi_{\pi}$ d'ordre fini.

\bigskip

On doit aussi introduire les paires $(Q,\pi)$ constitu\'ees d'un sous-groupe
paraboli\-que standard $Q \varb \GL_r$ associ\'e \`a une partition $r = r_1 +
\cdots + r_k$, avec pour sous-groupe de L\'evi $M_Q = \GL_{r_1} \ts \cdots \ts
\GL_{r_k}$, et d'une repr\'esentation irr\'eductible $\pi$ de $M_Q (\Ab)$ qui est
le produit $\pi_1 \ot \ldots \ot \pi_k$ de repr\'esentations automorphes
cuspidales $\pi_1 , \ldots , \pi_k$ de $\GL_{r_1} (\Ab) , \ldots , \GL_{r_k}
(\Ab)$.

Deux paires $(Q,\pi)$ et $(Q' , \pi')$ sont dites \'equivalentes s'il existe une
permutation $\s$ \'echangeant les facteurs de $M_Q = \GL_{r_1} \ts \ldots \ts
\GL_{r_k}$ et de $M_{Q'}$ et un caract\`ere $\lb$ de $M_Q (\Ab)$ ``non ramifi\'e''
c'est-\`a-dire se factorisant \`a travers l'homomorphisme
$$
\GL_{r_1} (\Ab) \ts \ldots \ts \GL_{r_k} (\Ab) \build\longra_{}^{\det}
(\Ab^{\ts})^k \build\longra_{}^{\deg} \Zb^k
$$
tels que $\pi' \cong \s (\lb \ot \pi)$. On choisit un repr\'esentant $(Q,\pi)$
dans chaque classe d'\'equivalence.

Soit $a \in \Ab^{\ts}$ un ad\`ele inversible de degr\'e non nul. L'espace de
Hilbert $L^2 (\GL_r (F) \bsh \GL_r (\Ab) / a^{\Zb})$ des fonctions de carr\'e
int\'egrable sur $\GL_r (F) \bsh \GL_r (\Ab) / a^{\Zb}$ est muni d'une action de
l'alg\`ebre de Hecke $\Hc^r$ par convolution \`a droite. Le th\'eor\`eme
fondamental de la th\'eorie des fonctions automorphes d\'ecrit sa
d\'ecomposi\-tion comme somme de repr\'esentations irr\'eductibles~:

\smallskip

{\bf Th\'eor\`eme de d\'ecomposition spectrale de Langlands.} -- \it  On a
\begin{eqnarray}
L^2 (\GL_r (F) \bsh \GL_r (\Ab) / a^{\Zb}) &= &\build\bigoplus_{{\pi \, \in \, \{ \pi
\}_r \atop \chi_{\pi} (a) \, = \, 1}}^{} \pi \nonumber \\
&\op &\build\bigoplus_{(Q,\pi)}^{} \ \hbox{$[s$omme continue de repr\'esentations}
\nonumber \\
&&\hbox{irr\'eductibles construites \`a partir de $(Q,\pi)$} \nonumber \\
&&\hbox{via la th\'eorie des s\'eries d'Eisenstein$]$}. \nonumber
\end{eqnarray}
\rm

\smallskip

Nous \'enon\c cons ici ce th\'eor\`eme en termes tr\`es vagues. En effet, on a
seulement besoin de savoir pour la suite que les repr\'esentations automorphes
cuspidales de $\GL_r$ apparaissent dans la somme et que tout le reste provient des
rangs $< r$.

\subsection{La correspondance de Langlands sur les corps de
fonctions}\label{ssec1d}\setzero \vskip-5mm \hspace{5mm }

On choisit un isomorphisme alg\'ebrique entre $\ov{\Qb}_{\ell}$ et $\Cb$.

Avec Langlands, disons qu'une repr\'esentation automorphe cuspidale $\pi$ de
$\GL_r$ et une repr\'esentation $\ell$-adique $\s$ de $G_F$ irr\'eductible de
dimension $r$ se correspondent si, en tout $x \in \vert X \vert$ o\`u $\pi$ est
non ramifi\'ee, $\s$ est non ramifi\'ee et
$$
{\rm L}_x (\s , T) = {\rm L}_x (\pi , T) \, .
$$

\smallskip

{\bf Th\'eor\`eme.} -- \it  Pour tout entier $r \geq 1$, cette correspondance
d\'efinit une bijection
$$
\{ \pi \}_r \cong \{ \s \}_r \qquad \pi \mpo \s_{\pi} \ , \qquad \s \mpo \pi_{\s}
\, .
$$
\rm

\smallskip

Le cas $r=1$ est la th\'eorie du corps de classes global pour $F = F(X)$.

Le cas $r=2$ a \'et\'e d\'emontr\'e par Drinfeld au d\'ebut des ann\'ees 70 en
inventant les chtoucas et \'etudiant ceux de rang 2.

Le cas $r \geq 3$ a \'et\'e d\'emontr\'e par l'auteur il y a deux ans en \'etudiant
les chtoucas de Drinfeld de rang $r$.

\smallskip

On sait que l'isomorphisme de Satake s'\'etend en une ``correspondance de Lan\-glands
locale'' qui, pour tous les corps $F_x$ localis\'es de $F=F(X)$ en les points $x \in
\vert X \vert$, a \'et\'e d\'emontr\'ee il y a dix ans par Laumon, Rapoport et Stuhler
en \'etudiant la cohomologie $\ell$-adique des vari\'et\'es modulaires de ``${\cal
D}$-faisceaux elliptiques''. Elle fait se correspondre les repr\'esentations
$\ell$-adiques de dimension $r$ du groupe de Galois $G_{F_x}$ de $F_x$ et les
repr\'esentations lisses admissibles irr\'eductibles de $\GL_r (F_x)$.

La correspondance globale du th\'eor\`eme ci-dessus est compatible avec la
correspondance locale au sens que si $\pi \in \{ \pi \}_r$ et $\s \in \{ \s \}_r$ se
correspondent au sens global, alors en tout point $x \in \vert X \vert$ (y compris ceux
o\`u il y a ramification), le facteur local $\pi_x$ de $\pi$ en $x$ et la restriction
$\s_x$ de $\s$ \`a $G_{F_x}$ se correspondent au sens local.

\section{Formule des traces d'Arthur-Selberg et chtoucas de Drinfeld}
\label{section 2} \setzero\vskip-5mm \hspace{5mm }

On va commencer \`a expliquer la d\'emonstration de la correspondance de Lan\-glands
en rang $r$. La premi\`ere chose importante \`a dire est qu'elle se fait par
r\'ecurrence. On suppose $r \geq 2$ et la correspondance d\'ej\`a connue en tous les
rangs $r' < r$.

\subsection{Des repr\'esentations galoisiennes vers les repr\'esentations
automorphes}\label{ssec2a}\setzero \vskip-5mm \hspace{5mm }

Pour $\s \in \{ \s \}_r$, l'unicit\'e de $\pi_{\s} \in \{ \pi \}_r$ lui correspondant
au sens de Langlands r\'esulte du ``th\'eor\`eme de multiplicit\'e 1 fort'' de
Piatetski-Shapiro qui dit qu'une repr\'esentation automorphe cuspidale de $\GL_r$ est
caract\'eris\'ee par la connaissance de ses valeurs propres de Hecke en presque tout $x
\in \vert X
\vert$.

Quant \`a l'existence de $\pi_{\s}$, il est connu depuis les ann\'ees 80 (avec la
``formule du produit'' de Laumon) qu'elle r\'esulte de la correspondance de Langlands
en les rangs $r' < r$. Rappelons comment~:

La th\'eorie des mod\`eles de Whittaker permet d'associer \`a $\s$ une
repr\'esentation lisse admissible irr\'eductible $\pi$ de $\GL_r$ qui, en tout $x$
o\`u $\s$ est non ramifi\'ee, est non ramifi\'ee et v\'erifie ${\rm L}_x (\pi , T) =
{\rm L}_x (\s , T)$, et qui est r\'ealis\'ee dans un espace de fonctions $\vp$ sur $P_1
(F) \bsh \GL_r (\Ab)$, o\`u $P_1 \varb \GL_r$ d\'esigne le sous-groupe parabolique
standard de type $r = (r-1) + 1$. Le probl\`eme est de montrer que toutes ces
fonctions $\vp$ sont invariantes \`a gauche par $\GL_r (F)$ tout entier.

D'apr\`es les ``th\'eor\`emes r\'eciproques'' de Hecke, Weil (pour $\GL_2$) et
Piatetski-Shapiro (pour $\GL_r$, $r \geq 3$), c'est \'equivalent \`a montrer que pour
tout $r' < r$ et toute repr\'esentation automorphe cuspidale $\pi' \in \{ \pi
\}_{r'}$, la fonction L globale
$$
{\rm L} (\s \ts \pi' , T)
$$
est un polyn\^ome qui v\'erifie une certaine \'equation fonctionnelle.

Or, d'apr\`es l'hypoth\`ese de r\'ecurrence, $\pi'$ correspond \`a une
repr\'esentation galoisienne $\s' \in \{ \s \}_{r'}$ et on peut \'ecrire
$$
{\rm L} (\s \ts \pi' , T) = {\rm L} (\s \ot \s' , T) \, .
$$
Comme on est sur un corps de fonctions $F=F(X)$, on sait gr\^ace \`a Grothendieck
qu'une telle fonction L galoisienne ${\rm L} (\s \ot \s' , T)$ est un polyn\^ome et
v\'erifie une \'equation fonctionnelle induite par la dualit\'e de Poincar\'e sur la
courbe $X$. Il faut encore v\'erifier que la constante dans l'\'equation
fonctionnelle est celle qu'on veut et ceci est la formule du produit de Laumon.

\subsection{Quotients ad\'eliques et fibr\'es vectoriels sur la
courbe}\label{ssec2b}\setzero \vskip-5mm \hspace{5mm }

Les diff\'erentes approches g\'eom\'etriques de Drinfeld pour le programme de
Langlands sur les corps de fonctions sont fond\'ees sur la remarque suivante
d'Andr\'e Weil~:

\smallskip

{\bf Lemme.} -- \it  Le quotient double $\GL_r (F) \bsh \GL_r (\Ab) / K$ s'identifie
\`a l'ensemble des classes d'isomorphie de fibr\'es vectoriels de rang $r$ sur la
courbe $X$.

Plus g\'en\'eralement, si $K_N \sbsq K = \GL_r (O_{\Ab})$ est le sous-groupe ouvert
associ\'e \`a un niveau $N \hra X$, $\GL_r (F) \bsh \GL_r (\Ab) / K_N$ s'identifie
\`a l'ensemble des classes de fibr\'es $\Ec$ de rang $r$ sur $X$ munis d'une
structure de niveau $N$ c'est-\`a-dire d'un isomorphisme $\Ec \ot_{\Oc_X} \Oc_N =
\Ec_N \build\longra_{}^{\sim} \Oc_X^r$. \rm

\smallskip

Autrement dit, $\GL_r (F) \bsh \GL_r (\Ab) / K$ est l'ensemble des points \`a valeurs
dans $\Fb_q$ du champ $\Vc ec^r$ des fibr\'es vectoriels de rang $r$ sur $X$.

Cette identification a permis \`a Drinfeld de construire l'application $\s \mpo
\pi_{\s}$ en rang $r=2$ d'une fa\c con purement g\'eom\'etrique, construction qui,
dans le cas non ramifi\'e $(N = \es)$, a \'et\'e progressivement g\'en\'eralis\'ee en
rang $r$ arbitraire par Laumon, Kazhdan, Frenkel, Gaitsgory, Vilonen.

On part d'une repr\'esentation $\ell$-adique partout non ramifi\'ee $\s \in \{ \s
\}_r$ qu'on voit comme un syst\`eme local sur $X$ suppos\'e absolument
irr\'eductible. En r\'einterpr\'e\-tant g\'eom\'etriquement la construction des
mod\`eles de Whittaker, on lui associe un certain complexe $\Aut'_{\s}$ de faisceaux
$\ell$-adiques sur le champ ${\Vc ec'}^r$ des fibr\'es $\Ec$ de rang $r$ sur $X$
munis d'un plongement $\Om_X^1 \hra \Ec$ du fibr\'e inversible canonique de $X$. Le
quotient $P_1 (F) \bsh \GL_r (\Ab) / K$ s'identifie \`a l'ensemble des points \`a
valeurs dans $\Fb_q$ d'un certain ouvert de ${\Vc ec'}^r$ et la fonction $\vp$ qui
associe \`a tout tel point la somme altern\'ee des traces de l'\'el\'ement de
Frobenius agissant sur la fibre de $\Aut'_{\s}$ en ce point est celle associ\'ee \`a
$\s$ par la th\'eorie classique des mod\`eles de Whittaker. Le probl\`eme est de
montrer que $\Aut'_{\s}$ ``se descend'' par le morphisme ${\Vc ec'}^r \ra \Vc ec^r$
d'oubli des plongements de $\Om_X^1$ c'est-\`a-dire qu'il est l'image r\'eciproque
d'un certain faisceau pervers $\Aut_{\s}$ sur $\Vc ec^r$. Or il y a un grand ouvert
de $\Vc ec^r$ au-dessus duquel le quotient de ${\Vc ec'}^r$ par $\Gb_m$ est un
fibr\'e projectif. Comme les espaces projectifs sont simplement connexes, la descente
se fait automatiquement si l'on sait que la restriction de $\Aut'_{\s}$ au-dessus de
cet ouvert est un syst\`eme local ($\Gb_m$-\'equivariant). Frenkel, Gaitsgory et
Vilonen ont montr\'e que cette propri\'et\'e r\'esulte de l'annulation de certains
foncteurs cohomologiques dits ``de moyennisation'' et Gaitsgory a r\'ecemment
d\'emontr\'e cet \'enonc\'e d'annulation par des arguments purement g\'eom\'etriques.

On voit que cette d\'emonstration de l'existence de l'application $\s \mpo \pi_{\s}$ est
profond\'ement
diff\'erente de celle du paragraphe pr\'ec\'edent car ici la descente repose sur la
simple connexit\'e des espaces projectifs et non plus sur les propri\'et\'es des
fonctions L d\'eduites de l'existence des applications $\pi' \mpo \s_{\pi'}$ en rangs
$<r$.

\subsection{La formule des traces d'Arthur-Selberg}\label{ssec2c}\setzero \vskip-5mm
\hspace{5mm }

On veut aborder maintenant la construction de l'application $\pi \mpo \s_{\pi}$.
L'unicit\'e des $\s_{\pi}$ correspondant aux $\pi \in \{ \pi \}_r$ r\'esulte du
th\'eor\`eme de densit\'e de Chebotarev.

Pour l'existence, il faut commencer par avoir une prise sur l'ensemble $\{ \pi \}_r$.
Une telle prise est fournie par la formule des traces d'Arthur-Selberg que nous
rappelons~:

L'image $\vp * h$ d'une fonction $\vp \in L^2 (\GL_r (F) \bsh \GL_r (\Ab) / a^{\Zb})$
par convolution par $h \in \Hc^r$ est donn\'ee par
$$
(\vp * h)(g') = \int_{\GL_r (F) \bsh \GL_r (\Ab) / a^{\Zb}} K_{h,G} (g',g) \, \vp (g)
\cdot {\rm d} g
$$
o\`u
$$
K_{h,G} (g',g) = \sum_{{\g \, \in \, \GL_r (F) \atop n \, \in \, \Zb}} h(g^{-1} \, \g
\, a^n \, g') \, .
$$
Pour $r \geq 2$, le quotient $\GL_r (F) \bsh \GL_r (\Ab) / a^{\Zb}$ a certes un
volume fini mais il n'est pas compact si bien que l'action de $h \in \Hc^r$ n'a pas
de trace au sens que l'int\'egrale
$$
\hbox{``}\Tr (h)\hbox{''} = \int_{\GL_r (F) \bsh \GL_r (\Ab) / a^{\Zb}} K_{h,G} (g,g)
\cdot {\rm d} g
$$
diverge en g\'en\'eral. Il r\'esulte du th\'eor\`eme de d\'ecomposition spectrale de
Langlands que le noyau $K_{h,G}$ s'\'ecrit naturellement
$$
K_{h,G} = \sum_{{\pi \, \in \, \{ \pi \}_r \atop \chi_{\pi} (a) \, = \, 1}}
K_{h,G}^{\pi} + \sum_{(Q,\pi)} K_{h,G}^{Q,\pi}
$$
et il faut pr\'eciser que le probl\`eme de divergence provient des $Q \varb G =
\GL_r$.

\smallskip

Pour cette raison, on introduit les troncatures d'Arthur~:

Pour tout sous-groupe parabolique standard $P \varb G$, l'action de $h$ sur \break
$L^2 (P(F) \bsh \GL_r (\Ab) / a^{\Zb})$ a un noyau
$$
(g',g) \mpo K_{h,P} (g',g) = \sum_{{\g \, \in \, P(F) \atop n \, \in \, \Zb}}
h(g^{-1} \, \g \, a^n \, g')
$$
qui se d\'ecompose naturellement en une somme
$$
K_{h,P} = \sum_{(Q,\pi)} K_{h,P}^{Q,\pi}
$$
o\`u n'apparaissent que les $Q$ tels que $M_Q \sbsq M_P$ \`a permutation pr\`es des
facteurs.

On consid\`ere un polygone de troncature, c'est-\`a-dire une fonction convexe $p :
[0,r] \ra \Rb_+$, s'annulant en 0 et $r$ et affine sur chaque intervalle $[r' - 1 ,
r']$, $0 < r' \leq r$. On pose
$$
K_{h,G}^{\leq p} (g,g) = K_{h,G} (g,g) + \sum_{P \, \varb \, \GL_r \, = \, G}
(-1)^{\vert P \vert - 1} \sum_{\d \, \in \, P(F) \bsh G(F)} \un_P^p (\d g) \, K_{h,P}
(\d g , \d g)
$$
o\`u $\vert P \vert$ d\'esigne le nombre de facteurs de chaque $M_P$ et $\un_P^p$ est
une fonction ca\-ract\'eristique (prenant les valeurs 1 ou 0) sur $P(F) \bsh G(\Ab)$
qui d\'epend du polygone de troncature $p$. Puis on pose
$$
\Tr^{\leq p} (h) = \int_{G(F) \bsh G(\Ab) / a^{\Zb}} K_{h,G}^{\leq p} (g,g) \cdot {\rm
d} g \, .
$$
On d\'efinit aussi des $\Tr_{Q,\pi}^{\leq p} (h)$ par des sommes altern\'ees et
int\'egrales semblables \`a partir des $K_{h,P}^{Q,\pi} (\cdot , \cdot)$.

Comme cons\'equence du th\'eor\`eme de d\'ecomposition spectrale de Langlands, on
montre~:

\smallskip

{\bf Th\'eor\`eme (formule des traces d'Arthur-Selberg).} -- \it  Toutes ces
int\'egrales convergent et on a
$$
\Tr^{\leq p} (h) = \sum_{{\pi \, \in \, \{ \pi \}_r \atop \chi_{\pi} (a) \, = \, 1}}
\Tr_{\pi} (h) + \sum_{(Q,\pi)} \Tr_{Q,\pi}^{\leq p} (h) \, .
$$
\rm

\smallskip

On voit donc que dans $\Tr^{\leq p} (h)$ apparaissent toutes les traces $\Tr_{\pi}
(h)$ plus d'autres termes (qui peuvent \^etre explicit\'es). Ceux-ci sont tr\`es
compliqu\'es, ils d\'ependent du polygone de troncature $p$ et ils ne sont pas
invariants par conjugaison de $h$ mais le plus important est qu'ils proviennent des
rangs $<r$.

Un autre point important est que les troncatures d'Arthur trouvent un sens
g\'eom\'etrique dans l'\'equivalence de Weil~:

\smallskip

{\bf Proposition.} -- \it  Si le polygone de troncature $p : [0,r] \ra \Rb_+$ est
``assez convexe'' $($c'est-\`a-dire si les diff\'erences de pentes $[p(r') - p(r'-1)]
- [p(r'+1) - p(r')]$, $0 < r' < r$, sont assez grandes$)$ en fonction de $h$, alors
pour tout $g \in G(F) \bsh G(\Ab)$ auquel est associ\'e un fibr\'e $\Ec$ de rang $r$
sur $X$, on a
$$
K_{h,G}^{\leq p} (g,g) = \Biggl\{ \begin{matrix}
K_{h,G} (g,g) &\hbox{si le polygone canonique de Harder-Narasimhan} \\
&\hbox{de $\Ec$ est $\leq p$}, \hfill \\
0 &\hbox{sinon}. \hfill
\end{matrix}
$$
\rm

\smallskip

\subsection{Les chtoucas de Drinfeld}\label{ssec2d}\setzero \vskip-5mm \hspace{5mm }

Pour \^etre exploit\'ee, la formule des traces d'Arthur-Selberg a besoin d'\^etre
combin\'ee avec autre chose. Les chtoucas de Drinfeld vont permettre de lui donner
une interpr\'etation g\'eom\'etrique et de l\`a une interpr\'etation cohomologique
via le th\'eor\`eme des points fixes de Grothendieck-Lefschetz.

Le probl\`eme est de construire des repr\'esentations $\ell$-adiques $\s_{\pi}$ de
$G_F$. En g\'eom\'etrie alg\'ebrique, on peut construire des repr\'esentations
galoisiennes en d\'efinis\-sant des vari\'et\'es sur $F$ (ou sur $X$ ou, comme il va
se produire, sur $X \ts X$) et en prenant leur cohomologie $\ell$-adique (ou celle de
leur fibre g\'en\'erique). Dans notre situa\-tion, il faut bien s\^ur que ces
vari\'et\'es aient un rapport \'etroit avec les repr\'esentations automorphes et donc
avec le quotient $\GL_r (F) \bsh \GL_r (\Ab)$.

Revenant \`a l'\'equivalence de Weil, on remarque que se donner un fibr\'e sur $X$
\'equivaut \`a se donner un fibr\'e $\Ec$ sur $\ov X = X \ot_{\Fb_q} \ov{\Fb}_q$ muni
d'un isomorphisme avec son transform\'e $^{\tau}\Ec = (\Id_X \ot \Frob)^* \Ec$ par
l'endomorphisme de Frobenius. Drinfeld a d\'ecouvert qu'en autorisant cet
isomorphisme $\Ec \cong \, ^{\tau}\Ec$ \`a avoir un p\^ole et un z\'ero, on d\'efinit
un probl\`eme de modules dont le classifiant r\'epond \`a la question pos\'ee~:

\smallskip

{\bf D\'efinition (Drinfeld).} -- (i) \it Un chtouca de rang $r$ \`a valeurs dans un
sch\'ema $S$ $($sur le corps de base $\Fb_q)$ consiste en

\smallskip

\noindent $\bu$ un fibr\'e $\Ec$ localement libre de rang $r$ sur $X \ts S$,

\smallskip

\noindent $\bu$ une modification de $\Ec$ c'est-\`a-dire un diagramme
$$
\Ec \build\lhra_{}^{j} \Ec' \build\lhla_{}^{t} \Ec''
$$
o\`u $\Ec' , \Ec''$ sont des fibr\'es sur $X \ts S$ et o\`u $j,t$ sont des
plongements dont les conoyaux sont support\'es par les graphes de deux morphismes
$\ify , 0 : S \ra X$ et sont inversibles comme faisceaux coh\'erents sur $\Oc_S$,

\smallskip

\noindent $\bu$ un isomorphisme $(\Id_X \ts \Frob_S)^* \Ec = \, ^{\tau}\Ec
\build\longra_{}^{\sim} \Ec''$. \rm

\smallskip

(ii) \it Une structure de niveau $N \hra X$ sur un tel chtouca $($dont le p\^ole
$\ify$ et le z\'ero $0$ \'evitent $N)$ est un isomorphisme $\Ec \ot_{\Oc_{X \ts S}}
\Oc_{N \ts S} = \Ec_N \build\longra_{}^{\sim} \Oc_{N \ts S}^r$ compatible avec
l'isomorphisme $^{\tau}\Ec_N \build\longra_{}^{\sim} \Ec_N$. \rm

\smallskip

Quand $S$ varie, les groupo{\"\i}des de chtoucas de rang $r$ [resp. et avec
structures de niveau $N$] constituent un champ $\Cht^r$ [resp. $\Cht_N^r$] dont on
montre qu'il est alg\'ebrique au sens de Deligne-Mumford et localement de type fini.
Voici les principales propri\'et\'es de ces champs modulaires~:

\smallskip

\noindent $\bu$ D'associer \`a tout chtouca son p\^ole et son z\'ero d\'efinit un
morphisme
$$
(\ify , 0) : \Cht^r \ra X \ts X \qquad \hbox{[resp.} \ \Cht_N^r \ra (X-N) \ts (X-N)]
$$
qui est lisse de dimension relative $2r-2$.

\smallskip

\noindent $\bu$ Chaque $\Cht_N^r$ est repr\'esentable fini \'etale galoisien de
groupe $\GL_r (\Oc_N)$ au-dessus de $\Cht^r \ts_{X \ts X} (X-N) \ts (X-N)$.

\smallskip

\noindent $\bu$ Chaque $\Cht_N^r$ est muni d'une action du groupe $F^{\ts} \bsh
\Ab^{\ts}$ et d'une action par correspondances finies \'etales de la sous-alg\`ebre
de Hecke $\Hc_N^r$.

\bigskip

On s'int\'eresse aux espaces de cohomologie $\ell$-adique \`a supports compacts
\break $H_c^{\nu} (\Cht_N^r / a^{\Zb})$, $0 \leq \nu \leq 2(2r-2)$, des $\Cht_N^r /
a^{\Zb}$ au-dessus du point g\'en\'erique $\Spec F^2$ de $X \ts X$. Ils sont munis
d'une double action du groupe de Galois $G_{F^2}$ et des alg\`ebres $\Hc_N^r$. Notre
but va \^etre de les calculer au moins partiellement et de montrer qu'ils r\'ealisent
l'application  $\pi \mpo \s_{\pi}$.

Chaque $\Cht_N^r / a^{\Zb}$ n'a qu'un nombre fini de composantes connexes (ce qui
correspond au fait que le volume de $\GL_r (F) \bsh \GL_r (\Ab) / a^{\Zb}$ est fini)
mais il n'est pas de type fini (de m\^eme que $\GL_r (F) \bsh \GL_r (\Ab) / a^{\Zb}$
n'est pas compact), ses espaces de cohomologie $\ell$-adique sont de dimension
infinie (de m\^eme que la d\'ecomposition spectrale de $L^2 (\GL_r (F) \bsh \GL_r (\Ab)
/ a^{\Zb}$) fait appara{\^\i}tre des sommes continues) et les correspondances de Hecke
compos\'ees avec les \'el\'ements de Frobenius ont une infinit\'e de points fixes (de
m\^eme que les int\'egrales ``$\Tr (h)$'' divergent).

Pour surmonter cette difficult\'e, on consid\`ere \`a nouveau un polygone de
troncature $p : [0,r] \ra \Rb_+$ et on demande aux chtoucas que leur polygone
canonique de Harder-Narasimhan (qui se d\'efinit facilement car un chtouca est un
fibr\'e muni d'une structure suppl\'ementaire) soit $\leq p$. Cela d\'efinit des
ouverts $\Cht_N^{r, \ov p \, \leq \, p}$ dans les $\Cht_N^r$ dont les quotients par
$a^{\Zb}$ sont de type fini (de m\^eme que les int\'egrales $\Tr^{\leq p} (h)$
convergent).

Une fonction $h \in \Hc_N^r$ contient en facteur la fonction caract\'eristique
$\un_y$ de $\GL_r (O_y)$ dans $\GL_r (F_y)$ en presque tout $y \in \vert X \vert$.
On consid\`ere un point ferm\'e $x$ de $(X-N) \ts (X-N)$ dont les deux projections
$\ify , 0 \in \vert X \vert$ sont distinctes et v\'erifient cette propri\'et\'e, et
un multiple $s = \deg (\ify) \, s' = \deg (0) \, u'$ de $\deg (x)$. Le compos\'e $h
\ts \Frob^s$ a dans la fibre de $\Cht_N^{r, \ov p \, \leq \, p} / a^{\Zb}$ au-dessus
de n'importe quel point de $(X \ts X)(\ov{\Fb}_q)$ support\'e par $x$ un nombre fini
de points fixes not\'e
$$
\Lef_x (h \ts \Frob^s , \Cht_N^{r, \ov p \, \leq \, p} / a^{\Zb}) \, .
$$
En utilisant une description ad\'elique des points fixes (rendue possible par la
pro\-ximit\'e des chtoucas avec le quotient $\GL_r (F) \bsh \GL_r (\Ab)$) et le calcul
des int\'egrales orbitales des fonctions sph\'eriques $h_{\ify}^{-s'}$ et $h_0^{u'}$
sur $\GL_r (F_{\ify})$ et $\GL_r (F_0)$ dont les transform\'es de Satake sont les
polyn\^omes sym\'etriques $q^{s \, \frac{(r-1)}{2}} (Z_1^{-s'} + \cdots + Z_r^{-s'})$
et $q^{s \, \frac{(r-1)}{2}} (Z_1^{u'} + \cdots + Z_r^{u'})$, on relie ce nombre \`a
la trace tronqu\'ee $\Tr^{\leq p} (h')$ de la fonction $h'$ d\'eduite de $h$ en
rempla\c cant les facteurs $\un_{\ify}$ et $\un_0$ par $h_{\ify}^{-s'}$ et
$h_0^{u'}$. Par combinaison avec la formule des traces d'Arthur-Selberg, on obtient~:

\smallskip

{\bf Th\'eor\`eme (comptage des points fixes).} -- \it Dans la situation ci-dessus et
si $p$ est assez convexe et $\deg (\ify)$, $\deg (0)$, $s$ sont assez grands en
fonction de $h$, on a
\begin{eqnarray}
&&\frac{1}{r!} \sum_{k \, = \, 1}^{r!} \Lef_{(\Frob_X^k \ts \Id_X)(x)} (h \ts
\Frob^s , \Cht_N^{r,\ov p \, \leq \, p} / a^{\Zb}) \nonumber \\
&= &q^{(r-1)s} \sum_{{\pi \, \in \, \{ \pi \}_r \atop \chi_{\pi} (a) \, = \, 1}}
\Tr_{\pi} (h) (z_1 (\pi_{\ify})^{-s'} + \cdots + z_r (\pi_{\ify})^{-s'}) \nonumber \\
&&\hspace{40mm} (z_1 (\pi_0)^{u'} + \cdots + z_r (\pi_0)^{u'}) \nonumber \\
&+ &\hbox{autres termes o\`u apparaissent les valeurs propres de Hecke en $\ify$ et
$0$ des} \nonumber \\
&&\hbox{repr\'esentations automorphes cuspidales des $\GL_{r_1} \ts \cdots \ts
\GL_{r_k}$}, \nonumber \\
&&r_1 + \cdots + r_k = r. \nonumber
\end{eqnarray}
\rm

\smallskip

Si dans cette formule il n'y avait que le terme principal et si on disposait sur
$\Cht_N^{r, \ov p \, \leq \, p} / a^{\Zb}$ d'un th\'eor\`eme des points fixes de
Grothendieck-Lefschetz interpr\'etant ces nombres comme la trace sur la cohomologie
d'une action des $h \ts \Frob_x^{-s / \deg (x)}$ (en notant $\Frob_x \in G_{F^2}$ un
\'el\'ement de Frobenius en $x$), on aurait une repr\'esentation de $\Hc_N^r \ts
G_{F^2}$ qui se d\'ecomposerait n\'ecessairement en
$$
\bigoplus_{{\pi \, \in \, \{ \pi \}_r \atop \chi_{\pi} (a) \, = \, 1}} (\pi \cdot
\un_N)
\ot [\s_{\pi} \boxtimes \check{\s}_{\pi}] \, (1-r)
$$
et on aurait construit l'application cherch\'ee $\pi \mpo \s_{\pi}$.

Mais le terme principal n'est pas seul, il y a aussi des termes compl\'ementaires qui
d\'ependent de $p$ et font appara{\^\i}tre les repr\'esentations automorphes
cuspidales en rangs $<r$, et surtout il n'y a m\^eme pas d'action de $\Hc_N^r$ car
les correspondances de Hecke de $\Cht_N^r / a^{\Zb}$ ne stabilisent pas les ouverts
$\Cht_N^{r, \ov p \, \leq \, p} / a^{\Zb}$ (ce qui correspond au fait que les traces
tronqu\'ees $\Tr^{\leq p} (h)$ ne sont pas invariantes par conjugaison).

\section{Donner un sens cohomologique au comptage des points fixes} \label{section 3}
\setzero\vskip-5mm \hspace{5mm }

On vient de dire que les ouverts $\Cht_N^{r, \ov p \, \leq \, p} / a^{\Zb}$ ne sont
pas stabilis\'es par les correspondances de Hecke si bien qu'en dehors du cas $h =
\un_N$ (c'est-\`a-dire de la seule action de $G_{F^2}$) les nombres $\Lef_x (h \ts
\Frob^s , \Cht_N^{r, \ov p \, \leq \, p} / a^{\Zb})$ n'ont pas a priori de sens
cohomologique. On va voir qu'en fait il est possible de leur en donner un plus
\'elabor\'e.

\subsection{Compactifications}\label{ssec3a}\setzero \vskip-5mm \hspace{5mm }

On compactifie les ouverts tronqu\'es $\Cht_N^{r, \ov p \, \leq \, p} / a^{\Zb}$ pour
retrouver des cor\-respondances de Hecke qui agissent sur la cohomologie. On commence
par le cas sans niveau~:

\smallskip

Pour tout polygone de troncature $p$ assez convexe en fonction du genre de la courbe
$X$, on construit un champ $\ov{\Cht^{r, \ov p \, \leq \, p}}$ alg\'ebrique au sens
d'Artin, dont les groupes d'automorphismes sont finis (mais ramifi\'es), qui est muni
d'une action de $F^{\ts} \bsh \Ab^{\ts}$ et d'un morphisme lisse sur $X \ts X$, qui
contient $\Cht^{r, \ov p \, \leq \, p}$ comme ouvert, dont le bord est un diviseur
\`a croisements normaux relatif et enfin dont le quotient par $a^{\Zb}$ est propre
(en particulier de type fini et s\'epar\'e) sur $X \ts X$. Les strates de bord
$\Cht_{\ud r}^{r, \ov p \, \leq \, p}$ sont naturellement index\'ees par les
partitions $\ud r = (r = r_1 + \cdots + r_k)$ de l'entier $r$. Si on distingue le
degr\'e $d$ des chtoucas avec donc $\Cht^r = \build\coprod_{d \, \in \, \Zb}^{}
\Cht^{r,d}$, $\ov{\Cht^{r, \ov p \, \leq \, p}} = \build\coprod_{d \, \in \, \Zb}^{}
\ov{\Cht^{r,d, \ov p \, \leq \, p}}$ et $\Cht_{\ud r}^{r, \ov p \, \leq \, p} =
\build\coprod_{d \, \in \, \Zb}^{} \Cht_{\ud r}^{r,d, \ov p \, \leq \, p}$, chaque
$\Cht_{\ud r}^{r,d, \ov p \, \leq \, p}$ est essentiellement de la forme
$$
\Cht^{r_1 , d_1 , \ov p \, \leq \, p_1} \ts_X \Cht^{r_2 , d_2 , \ov p \, \leq \, p_2}
\ts_{X,\Frob} \cdots \ts_{X,\Frob} \Cht^{r_k , d_k , \ov p \, \leq \, p_k}
$$
o\`u $d_1 , \ldots , d_k$ et $p_1 , \ldots , p_k$ sont des degr\'es et polygones de
troncatures qui se d\'eduisent de $d$ et $p$. La strate $\Cht_{\ud r}^{r, \ov p \,
\leq \, p}$ est munie d'un morphisme vers $X \ts X^{k-1} \ts X$ o\`u le premier
facteur $X$ est le p\^ole, le dernier facteur $X$ est le z\'ero et les $k-1$ facteurs
$X$ suppl\'ementaires sont appel\'es les d\'eg\'en\'erateurs.

Toute fonction $h \in \Hc_{\es}^r$ induit une correspondance sur $\Cht^r / a^{\Zb}$.
On peut consid\'erer sa trace dans $(\Cht^{r, \ov p \, \leq \, p} / a^{\Zb})^2$ puis
la normalisation de celle-ci sur $(\ov{\Cht^{r, \ov p \, \leq \, p}} /$ $a^{\Zb})^2$.
Comme $\ov{\Cht^{r, \ov p \, \leq \, p}} / a^{\Zb}$ est propre et lisse sur $X \ts
X$, elle agit sur sa cohomologie.

\smallskip

Consid\'erons maintenant un niveau $N \hra X$. La normalisation $\ov{\Cht_N^{r, \ov p
\, \leq \, p}} / a^{\Zb}$ de $\ov{\Cht^{r, \ov p \, \leq \, p}} / a^{\Zb} \ts_{X \ts
X} (X-N) \ts (X-N)$ dans $\Cht_N^{r, \ov p \, \leq \, p} / a^{\Zb}$ contient
$\Cht_N^{r, \ov p \, \leq \, p} / a^{\Zb}$ comme ouvert et elle est propre sur $(X-N)
\ts (X-N)$ mais elle n'est pas lisse et l'auteur ignore comment r\'esoudre ses
singularit\'es. Toutefois, en demandant que non seulement le p\^ole et le z\'ero mais
aussi les d\'eg\'en\'erateurs \'evitent $N$, on d\'efinit un ouvert $\ov{\Cht_N^{r,
\ov p \, \leq \, p}}' / a^{\Zb}$ de $\ov{\Cht_N^{r, \ov p \, \leq \, p}} / a^{\Zb}$
qui contient strictement $\Cht_N^{r, \ov p \, \leq \, p} / a^{\Zb}$, est lisse sur $X
\ts X$ et dont le bord est un diviseur \`a croisements normaux relatif. Les
correspondances de Hecke $h \in \Hc_N^r$ agissent sur sa cohomologie d'apr\`es le
th\'eor\`eme suivant~:

\smallskip

{\bf Th\'eor\`eme de stabilit\'e globale.} -- \it Les correspondances de Hecke
\'etendues par normalisation sur $(\ov{\Cht_N^{r, \ov p \, \leq \, p}}' / a^{\Zb})^2$
stabilisent $\ov{\Cht_N^{r, \ov p \, \leq \, p}}' / a^{\Zb}$ au sens que leurs deux
projections sur celui-ci sont propres. \rm

\smallskip

Ce r\'esultat correspond sans doute au ph\'enom\`ene suivant dans la formule des
traces d'Arthur-Selberg~: Le d\'efaut d'invariance par conjugaison des traces
tronqu\'ees $\Tr^{\leq p} (h)$ se mesure par les d\'eriv\'ees logarithmiques des
op\'erateurs d'entre\-lacement de Langlands. Or ces op\'erateurs sont des produits
eul\'eriens et leurs d\'eriv\'ees logarithmiques sont des sommes sur tous les points
$x \in \vert X \vert$ lesquels correspondent exactement aux valeurs possibles des
d\'eg\'en\'erateurs. Il n'y a pas d'entrelacement donc pas d'instabilit\'e d'un point
de $\vert X \vert$ \`a un autre et en demandant au p\^ole, au z\'ero et aux
d\'eg\'en\'erateurs d'\'eviter certains points de $\vert X \vert$, on garde la
propri\'et\'e de stabilit\'e globale que v\'erifie automatiquement la
compactification toute enti\`ere $\ov{\Cht_N^{r, \ov p \, \leq \, p}} / a^{\Zb}$.

\subsection{Cohomologie n\'egligeable et cohomologie
essentielle}\label{ssec3b}\setzero \vskip-5mm \hspace{5mm }

En r\'esum\'e, on a maintenant toutes les structures et informations qu'on \break
pourrait
souhaiter mais elles sont dispers\'ees entre les diff\'erents objets $\Cht_N^r /
a^{\Zb}$, \break $\Cht_N^{r, \ov p \, \leq \, p} / a^{\Zb}$ et $\ov{\Cht_N^{r, \ov p
\, \leq \, p}}' / a^{\Zb}$~:

Sur $\Cht_N^r / a^{\Zb}$, on a une action de l'alg\`ebre de Hecke $\Hc_N^r$ mais la
cohomologie est de dimension infinie et les ensembles de points fixes sont infinis.

L'ouvert $\Cht_N^{r, \ov p \, \leq \, p} / a^{\Zb}$ est de type fini et on y a une
formule de comptage des points fixes qui s'exprime en termes automorphes mais on a
perdu l'action des correspondances de Hecke.

Enfin, sur $\ov{\Cht_N^{r, \ov p \, \leq \, p}}' / a^{\Zb}$ et sa cohomologie il y a
\`a nouveau une action de chaque $h \in \Hc_N^r$ mais on ne donne pas de comptage des
points fixes et surtout il n'y a pas d'action de l'alg\`ebre $\Hc_N^r$ car la
normalisation des correspondances de Hecke ne commute pas avec la multiplication.

\bigskip

L'id\'ee est de d\'efinir dans toute repr\'esentation $\ell$-adique de $G_{F^2}$
(apr\`es semi-simplification) une partie ``n\'egligeable'' et une partie
``essentielle'', de fa\c con que les $\Cht_N^r / a^{\Zb}$, $\Cht_N^{r, \ov p \, \leq
\, p} / a^{\Zb}$ et $\ov{\Cht_N^{r, \ov p \, \leq \, p}}' / a^{\Zb}$ aient m\^eme
cohomologie essentielle (ce qui permettra de rassembler sur celle-ci les informations
dispers\'ees dont on dispose) et que les traces des actions sur la cohomologie
essentielle soient donn\'ees par les termes principaux dans la formule de comptage
des points fixes, ceux associ\'es aux repr\'esentations automorphes cuspidales $\pi
\in \{ \pi \}_r$ de $\GL_r$.

Or on s'attend \`a ce que les $\pi \in \{ \pi \}_r$ correspondent \`a des
repr\'esentations $\ell$-adiques de $G_F$ irr\'eductibles de dimension $r$. Quant aux
termes compl\'ementaires dans la formule de comptage, ils sont associ\'es aux
repr\'esentations automorphes cuspidales en rangs $<r$ lesquelles correspondent,
d'apr\`es l'hypoth\`ese de r\'ecurrence, aux repr\'esentations $\ell$-adiques de
$G_F$ irr\'eductibles de dimension $<r$. Cela dicte la d\'efinition suivante (o\`u
$q', q'' : G_{F^2} \rra G_F$ sont les deux homomorphismes induits par $X \ts X \rra
X$)~:

\smallskip

{\bf D\'efinition.} -- \it Une repr\'esentation $\ell$-adique irr\'eductible de
$G_{F^2}$ est dite ``$r$-n\'egligeable'' si elle est facteur direct d'une
repr\'esentation de la forme
$$
{q'}^* \, \s' \ot {q''}^* \, \s''
$$
avec $\s' , \s''$ deux repr\'esentations de $G_F$ irr\'eductibles de rangs $<r$.

Elle est dite ``essentielle'' sinon. \rm

\smallskip

\subsection{S\'eparation et identification de la cohomologie
essentielle}\label{ssec3c}\setzero \vskip-5mm \hspace{5mm }

Dans un premier temps, on oublie compl\`etement les correspondances de Hecke et on ne
consid\`ere que l'action de $G_{F^2}$ sur les diff\'erents espaces de cohomologie
$\ell$-adique \`a supports compacts. On montre~:

\smallskip

{\bf Proposition.} -- \it Pour tout niveau $N \hra X$, la cohomologie de $\Cht_N^r /
a^{\Zb}$ et des $\Cht_N^{r, \ov p \, \leq \, p} / a^{\Zb}$ et $\ov{\Cht_N^{r, \ov p
\, \leq \, p}}' / a^{\Zb}$ $($et la cohomologie d'intersection des $\ov{\Cht_N^{r,
\ov p \, \leq \, p}} / a^{\Zb}$ si $N \ne \es)$ ont la m\^eme partie essentielle
$H_N^{\ess}$. Elle est concentr\'ee en degr\'e m\'edian $\nu = 2r-2$ et pure de poids
$2r-2$. Si $x$ est un point ferm\'e de $(X-N) \ts (X-N)$ dont les deux projections
$\ify , 0 \in \vert X \vert$ sont distinctes et $s = \deg (\ify) \, s' = \deg (0) \,
u'$ est un multiple de $\deg (x)$, on a
\begin{eqnarray}
\Tr_{H_N^{\ess}} (\Frob_x^{-s / \deg (x)}) &= &q^{(r-1)s} \sum_{{\pi \, \in \, \{ \pi
\}_r \atop \chi_{\pi} (a) \, = \, 1}} \dim (\pi \cdot \un_N) (z_1 (\pi_{\ify})^{-s'} +
\cdots + z_r (\pi_{\ify})^{-s'}) \nonumber \\
&&\hspace{45mm} (z_1 (\pi_0)^{u'} + \cdots + z_r (\pi_0)^{u'}) \, . \nonumber
\end{eqnarray}
\rm

\smallskip

Quand $N = \es$, la cohomologie du bord $\ov{\Cht^{r, \ov p \, \leq \, p}} / a^{\Zb}
- \Cht^{r, \ov p \, \leq \, p} / a^{\Zb}$ est $r$-n\'egligeable car il est r\'eunion
de strates $\Cht_{\ud r}^{r, \ov p \, \leq \, p} / a^{\Zb}$ index\'ees par les
partitions $\ud r = (r=r_1 + \cdots + r_k)$, $k \geq 2$, qui se d\'evissent en termes
de $\Cht^{r_1} , \ldots , \Cht^{r_k}$. Pour $N \ne \es$, un argument plus
sophistiqu\'e utilisant ce d\'evissage montre aussi que les diff\'erences entre
$\Cht_N^{r, \ov p \, \leq \, p} / a^{\Zb}$, $\ov{\Cht_N^{r, \ov p \, \leq \, p}}' /
a^{\Zb}$ et $\ov{\Cht_N^{r, \ov p \, \leq \, p}} / a^{\Zb}$ sont $r$-n\'egligeables.

La formule de comptage des points fixes du paragraphe 2.4 (avec $h = \un_N$) et le
th\'eor\`eme des points fixes de Grothendieck-Lefschetz donnent les traces des
\'el\'ements $\Frob_x^{-s / \deg (x)}$ agissant sur la cohomologie des $\Cht_N^{r,
\ov p \, \leq \, p} / a^{\Zb}$. La s\'eparation de la partie essentielle se fait en
``testant'' ces espaces de cohomologie contre des repr\'esentations irr\'eductibles
$r$-n\'egligeables arbitraires et en regardant les p\^oles des fonctions L de
paires obtenues. Cela utilise la correspondance de Langlands d\'ej\`a connue en rangs
$<r$ par hypoth\`ese de r\'ecurrence et les propri\'et\'es classiques des fonctions
L de paires tant du c\^ot\'e automorphe que galoisien (en particulier
l'interpr\'etation cohomologique de Grothendieck et le th\'eor\`eme de puret\'e de
Deligne).

Enfin, $\Cht_N^r / a^{\Zb}$ et les $\Cht_N^{r, \ov p \, \leq \, p} / a^{\Zb}$ ont la
m\^eme cohomologie essentielle car la formule obtenue pour les traces des
$\Frob_x^{-s / \deg (x)}$ ne d\'epend pas de $p$.

\subsection{Action et traces des correspondances de Hecke}\label{ssec3d}\setzero
\vskip-5mm \hspace{5mm }

Dans un second temps, on revient aux correspondances de Hecke $h \in \Hc_N^r$.

On met d'abord une action naturelle de l'alg\`ebre $\Hc_N^r$ sur la cohomologie
essentielle $H_N^{\ess}$ en prouvant que $H = H_c^{2r-2} (\Cht_N^r / a^{\Zb}) = \
\build\varinjlim_{p}^{} H_c^{2r-2} (\Cht_N^{r, \ov p \, \leq \, p} / a^{\Zb})$ admet
une filtration finie $0 = H_0 \varb \ldots \varb H_i \varb \ldots \varb H_k = H$
respect\'ee par la double action de $G_{F^2}$ et de $\Hc_N^r$ dont les gradu\'es
impairs $H_{2i+1} / H_{2i}$ sont enti\`erement $r$-n\'egligeables et les gradu\'es
pairs $H_{2i+2} / H_{2i+1}$ sont enti\`erement essentiels.

\smallskip

Pour terminer, il reste \`a montrer que pour $h \in \Hc_N^r$, $x$ et $s = \deg (\ify)
\, s' = \deg (0) \, u'$ comme dans l'\'enonc\'e de la formule de comptage, on a
\begin{eqnarray}
\Tr_{H_N^{\ess}} (h \ts \Frob_x^{-s / \deg (x)}) &= &q^{(r-1)s} \sum_{{\pi \, \in \,
\{ \pi \}_r \atop \chi_{\pi} (a) \, = \, 1}} \Tr_{\pi} (h) (z_1 (\pi_{\ify})^{-s'} +
\cdots + z_r (\pi_{\ify})^{-s'}) \nonumber \\
&&\hspace{38mm} (z_1 (\pi_0)^{u'} + \cdots + z_r (\pi_0)^{u'}) \, . \nonumber
\end{eqnarray}
Pour cela, on a besoin d'une formule de Grothendieck-Lefschetz qui relie les nombres
de points fixes $\Lef_x (h \ts \Frob^s , \Cht_N^{r, \ov p \, \leq \, p} / a^{\Zb})$
\`a l'action de $h \ts \Frob^s$ sur la cohomologie de $\ov{\Cht_N^{r, \ov p \, \leq \,
p}}' / a^{\Zb}$ (et ce sera suffisant car on n'a plus alors qu'\`a combiner une telle
formule avec la formule de comptage et \`a identifier ce qui est ``essentiel'' des
deux c\^ot\'es par les arguments de fonctions L de paires).

On prouve en fait que pour toute $h \in \Hc_N^r$ fix\'ee, il existe des
correspondances cohomologiques ${\rm cl} (h)_{\ud r}$ agissant sur la cohomologie des
strates $\ov{\Cht_{N,\ud r}^{r, \ov p \, \leq \, p}}' / a^{\Zb}$ du bord
$\ov{\Cht_N^{r, \ov p \, \leq \, p}}' / a^{\Zb} - \Cht_N^{r, \ov p \, \leq \, p} /
a^{\Zb}$ telles que pour tout $x$ et tout $s$ comme plus haut on ait (en notant $H_c^*
(\cdot) = \, \build\sum_{\nu}^{} (-1)^{\nu} \, H_c^{\nu} (\cdot)^{\rm ss}$)~:
\begin{eqnarray}
&&\Lef_x (h \ts \Frob_s , \Cht_N^{r, \ov p \, \leq \, p} / a^{\Zb}) = \Tr_{H_c^*
(\ov{\Cht_N^{r, \ov p \, \leq \, p}}' / a^{\Zb})} (h \ts \Frob_x^{-s / \deg (x)})
\nonumber \\
&+ &\sum_{{\ud r \, = \, (r=r_1 + \cdots + r_k) \atop k \, \geq \, 2}} (-1)^{k-1}
\Tr_{H_c^* (\ov{\Cht_{N,\ud r}^{r, \ov p \, \leq \, p}}' / a^{\Zb})} ({\rm cl}
(h)_{\ud r} \ts \Frob_x^{-s / \deg (x)}) \, . \nonumber
\end{eqnarray}
La preuve de cette formule (qui justifie a posteriori le fait \'etrange qu'il y ait
une formule de comptage des points fixes dans l'ouvert non stable $\Cht_N^{r, \ov p \,
\leq \, p} / a^{\Zb}$) repose sur la propri\'et\'e g\'eom\'etrique suivante des
correspondances de Hecke (qui para{\^\i}t donc li\'ee \`a l'existence m\^eme de la
formule des traces d'Arthur-Selberg)~:

\smallskip

{\bf Th\'eor\`eme de ``stabilit\'e locale''.} -- \it  Dans $\hbox{{\goth X}} =
\ov{\Cht_{N}^{r, \ov p \, \leq \, p}}' / a^{\Zb}$, les correspondances de Hecke $h$
stabilisent l'ouvert $\hbox{{\goth X}}_{\es} = \Cht_{N}^{r, \ov p \, \leq \, p} /
a^{\Zb}$ ``localement au voisinage de leurs points fixes''.

Cela signifie que si $\G
\build{\lhook\joinrel\relbar\joinrel\relbar\joinrel\relbar\joinrel\rightarrow}_{}^{(p'
_{\G} , p''_{\G})} \hbox{{\goth X}} \ts \hbox{{\goth X}}$ est le
cycle qui supporte $h$, il existe un ouvert $U \sbs \hbox{{\goth
X}} \ts \hbox{{\goth X}}$ contenant les points fixes de toutes les
correspondances $\G \ts \Frob^n$, $n \in \Nb$, tel que
$$
{p'}_{\G}^{-1} (\hbox{{\goth X}}_{\es}) \cap U \sbsq
{p''}_{\G}^{-1} (\hbox{{\goth X}}_{\es}) \cap U \, .
$$
\rm

\section{Applications du th\'eor\`eme} \label{section 4}
\setzero\vskip-5mm \hspace{5mm }

La correspondance de Langlands entre repr\'esentations $\ell$-adiques irr\'eductibles
du groupe de Galois $G_F$ de $F$ et repr\'esentations automorphes cuspidales des
groupes $\GL_r$ sur $F$ a des cons\'equences imm\'ediates et importantes dans les deux
sens. Voici les principales~:

\subsection{\boldmath Cons\'equences sur les faisceaux $\ell$-adiques}\label{ssec4a}\setzero
\vskip-5mm \hspace{5mm }

On a d'abord des cons\'equences tr\`es fortes dans le cas des courbes~:

\smallskip

{\bf Th\'eor\`eme.} -- \it  Soient $X'$ un ouvert de la courbe $X$ et $\s$ un faisceau
$\ell$-adique lisse sur $X'$, qui est irr\'eductible de rang $r$ et dont le
d\'eterminant est un caract\`ere d'ordre fini. Alors~: \rm

\smallskip

\noindent (i) \it Il existe un corps de nombres $E \sbs \ov{\Qb}$ tel qu'en tout point
ferm\'e $x \in \vert X' \vert$, le polyn\^ome ${\rm L}_x (\s , T)^{-1} = \det_{\s} (1 - T
\cdot \Frob_x^{-1})$ soit \`a coefficients dans $E$. \rm

\smallskip

\noindent (ii) \it En tout $x \in \vert X' \vert$, les racines du polyn\^ome $\det_{\s}
(1-T \cdot \Frob_x^{-1})$ sont des nombres alg\'ebriques dont toutes les images
complexes sont de module $1$. Ce sont des unit\'es $\lb$-adiques en toutes les places
$\lb$ non archim\'ediennes et premi\`eres \`a $q$ de $E$. \rm

\smallskip

\noindent (iii) \it Pour toute place $\lb$ de $E$ au-dessus d'un nombre premier $\ell'$
ne divisant pas $q$, il existe sur $X'$ un faisceau $\ell'$-adique $\s_{\lb}$ lisse et
irr\'eductible de rang $r$ tel que
$$
{\textstyle\det_{\s}} (1 - T \cdot \Frob_x^{-1}) = {\textstyle\det_{\s_{\lb}}} (1 - T
\cdot \Frob_x^{-1}) \, , \qquad \fl x \in \vert X' \vert \, .
$$
\rm

\smallskip

On voit qu'en dimension 1 sur $\Fb_q$, un faisceau $\ell$-adique irr\'eductible dont le
d\'eterminant est d'ordre fini est ``pur de poids 0'' et ``il ne d\'epend pas du choix
de $\ell$''. La premi\`ere de ces deux propri\'et\'es s'\'etend automatiquement en
dimension arbitraire~:

\smallskip

{\bf Corollaire.} -- \it  Soient $\hbox{{\goth X}}$ une vari\'et\'e normale de type
fini sur $\Fb_q$ et $\s$ un faisceau $\ell$-adique lisse sur $\hbox{{\goth X}}$ qui est
irr\'eductible et dont le d\'eterminant est d'ordre fini.

Alors en tout point ferm\'e $x$ de $\hbox{{\goth X}}$, les valeurs propres de Frobenius
de $\s$ sont des nombres alg\'ebriques dont toutes les images complexes sont de modules
$1$.

Et ce sont des unit\'es $\lb$-adiques pour toute place $\lb$ premi\`ere \`a $q$. \rm

\subsection{Cons\'equences sur les repr\'esentations automorphes}\label{ssec4b}\setzero
\vskip-5mm \hspace{5mm }

Pour les repr\'esentations automorphes cuspidales des groupes lin\'eaires sur $F =
F(X)$, on a d'abord des cons\'equences sur les valeurs propres de Hecke et les
fonctions ${\rm L}$ de paires~:

\smallskip

{\bf Th\'eor\`eme.} -- (i) (Conjecture de Ramanujan-Petersson) \it Pour toute $\pi \in
\{ \pi \}_r$, les facteurs locaux $\pi_x$ de $\pi$ en les $x \in \vert X \vert$ sont
temp\'er\'es.

En particulier, en les $x$ o\`u $\pi_x$ est non ramifi\'e, ses valeurs propres de Hecke
v\'erifient
$$
\vert z_i (\pi_x) \vert = 1 \, , \qquad 1 \leq i \leq r \, .
$$
\rm

\noindent (ii) (Hypoth\`ese de Riemann g\'en\'eralis\'ee) \it Pour toute paire $\pi \in
\{ \pi \}_r$, $\pi' \in \{ \pi \}_{r'}$, tous les z\'eros de la fonction ${\rm L}$
globale
$$
{\rm L} (\pi \ts \pi' , T)
$$
sont sur le cercle
$$
\vert T \vert = q^{-1/2} \, .
$$
\rm

\smallskip

La partie (ii) du th\'eor\`eme est la traduction en termes automorphes du th\'eor\`eme de
puret\'e de Deligne.

D'autre part, on a les cas particuliers suivants de la fonctorialit\'e de Langlands
(o\`u, pour toute extension finie $F'$ de $F$, $\Ab_{F'}$ d\'esigne son anneau des
ad\`eles et $X_{F'}$ la courbe projective lisse associ\'ee qui est un rev\^etement fini
de $X = X_F$)~:

\smallskip

{\bf Th\'eor\`eme.} -- (i) (Existence du produit tensoriel automorphe) \it Soient $\pi$
et $\pi'$ deux repr\'esentations automorphes cuspidales de $\GL_r (\Ab_F)$ et $\GL_{r'}
(\Ab_F)$. Alors il existe une partition $rr' = r_1 + \cdots + r_k$ et des
repr\'esentations automorphes cuspidales $\pi^1 , \ldots , \pi^k$ de $\GL_{r_1} (\Ab_F)
, \ldots , \GL_{r_k} (\Ab_F)$ qui, en tout $x \in \vert X_F \vert$ o\`u $\pi$ et $\pi'$
sont non ramifi\'ees, sont elles-m\^emes non ramifi\'ees et v\'erifient
$$
\{ z_j (\pi_x) \, z_{j'} (\pi'_x) \mid 1 \leq j \leq r \, , \ 1 \leq j' \leq r' \} =
\coprod_{1 \, \leq \, i \, \leq \, k} \{ z_1 (\pi_x^i) , \ldots , z_{r_i} (\pi_x^i) \}
$$
\rm

\noindent (ii) (Changement de base) \it Soient $F'$ une extension finie de $F$ et $\pi$
une repr\'esentation automorphe cuspidale de $\GL_r (\Ab_F)$. Alors il existe une
partition $r = r_1 + \cdots + r_k$ et des repr\'esentations automorphes cuspidales
${\pi'}^1 , \ldots , {\pi'}^k$ de $\GL_{r_1} (\Ab_{F'}) , \ldots , \GL_{r_k}$
$(\Ab_{F'})$, non ramifi\'ees en tout point $x' \in \vert X_{F'} \vert$ de degr\'e
$\frac{\deg (x')}{\deg (x)}$ au-dessus d'un point $x \in \vert X_F \vert$ o\`u $\pi$
est non ramifi\'ee et qui v\'erifient
$$
\bigl\{ z_1 (\pi_x)^{\frac{\deg (x')}{\deg (x)}} , \ldots , z_r (\pi_x)^{\frac{\deg
(x')}{\deg (x)}} \bigl\} = \coprod_{1 \, \leq \, i \, \leq \, k} \{ z_1 ({\pi'}_x^i) ,
\ldots , z_{r_i} ({\pi'}_x^i) \} \, .
$$
\rm

\noindent (iii) (Induction automorphe) \it  Soient $F'$ une extension de $F$ de degr\'e
$d$ et $\pi'$ une repr\'esentation automorphe cuspidale de $\GL_r (\Ab_{F'})$. Alors il
existe une partition $rd = r_1 + \cdots + r_k$ et des repr\'esentations automorphes
cuspidales $\pi^1 , \ldots , \pi^k$ de $\GL_{r_1} (\Ab_F) , \ldots , \GL_{r_k}
(\Ab_F)$, non ramifi\'ees en tout point $x \in \vert X_F \vert$ au-dessus duquel les
points $x' \in \vert X_{F'} \vert$ sont non ramifi\'ees sur $x$ et pour $\pi'$ et qui
v\'erifient
$$
\prod_{1 \, \leq \, i \, \leq \, k} {\rm L}_x (\pi^i , T) = \prod_{x' \mid x} {\rm
L}_{x'} (\pi' , T) \, .
$$
\rm

\label{lastpage}

\end{document}